\documentclass[11pt,a4paper]{article}
\usepackage[english]{babel}
\usepackage{hyperref}
\usepackage[utf8]{inputenc}

\usepackage{color}

\usepackage{fontenc}

\usepackage{amssymb, amsmath, amsbsy}
\usepackage{array}
\usepackage{amsthm}

\usepackage{fancyhdr}

\usepackage{nccmath}

\usepackage{geometry}

\usepackage{graphicx} 
\usepackage{subfigure} 
\usepackage{float} 

\newcommand{\com}[1]{\opt{draft}{\textcolor{red}{
$\LHD$ #1 $\RHD$\marginpar{\textcolor{red}{$\begin{lema}acksquare$}}}}}

\newcommand{\comb}[1]{\opt{draft}{\textcolor{blue}{
$\LHD$ #1 $\RHD$\marginpar{\textcolor{blue}{$\begin{lema}acksquare$}}}}}

\newenvironment{demo}{{\bf Proof }}
{\qed \\}

\newcommand{\re}{\mathbb R}

\newcommand{\es}{\mathbb S}

\newcommand{\lb}{\label}
\newcommand{\nn}{\nonumber}

\newcommand{\eps}{\ensuremath{\varepsilon}}
\def\a{\alpha}

\def\dist{{\rm dist}}

\newcommand{\bde}{\begin{defi}}
\newcommand{\ede}{\end{defi}}

\newcommand{\be}{\begin{enumerate}}
\newcommand{\ee}{\end{enumerate}}

\newcommand{\ba}{\begin{array}}
\newcommand{\ea}{\end{array}}

\def\H{{\mathcal H}}

\def\dist{{\rm dist}}

\newtheorem{defi}{\hspace{12pt} Definition}
\newtheorem{teor}{\hspace{12pt} Theorem}
\newtheorem{prop}[teor]{\hspace{12pt} Proposition}

\newtheorem{lema}{\hspace{12pt} Lemma}
\newtheorem*{lema*}{\hspace{12pt} Lemma}
\newtheorem*{teor*}{\hspace{12pt} Theorem}

\newtheorem*{nota*}{\hspace{12pt} Remark}

\def\dist{{\rm dist}}

\def\H{{\mathcal H}}

\newcommand{\ben}{\begin{enumerate}}
\newcommand{\een}{\end{enumerate}}
\newcommand{\bi}{\begin{itemize}}
\newcommand{\ei}{\end{itemize}}
\newcommand{\bec}{\begin{equation}}
\newcommand{\eec}{\end{equation}}
\newcommand{\beca}{\begin{equation*}}
\newcommand{\eeca}{\end{equation*}}
\newcommand{\bal}{\begin{align}}
\newcommand{\aal}{\end{align}}
\newcommand{\bala}{\begin{align*}}
\newcommand{\aala}{\end{align*}}

\begin{document}

\title{Geodesic loops on tetrahedra in spaces of constant sectional curvature \footnotetext{The research is partially supported by grant PID2019-105019GB-C21 funded by MCIN/AEI/ 10.13039/501100011033 and by \lq\lq ERDF A way of makimg Europe\rq\rq,, and by the grant AICO 2021 21/378.01/1 funded by the Generalitat Valenciana.}
}
\author{Alexander Borisenko and Vicente Miquel}

\date{}

\maketitle

\begin{abstract}
Geodesic loops on polyhedra were studied only for Euclidean space and it was known that there are no simple geodesic loops on  regular tetrahedra. Here we prove that: 1) On the spherical space, there are no simple geodesic loops on tetrahedra with internal angles $\pi/3 < \a_i<\pi/2$  or regular tetrahedra with $\a_i=\pi/2$, and there are three simple geodesic loops for each vertex of a tetrahedra  with $\a_i > \pi/2$ and the lengths of the edges  $a_i>\pi/2$. 
2) On the hyperbolic space, for every regular tetrahedron $T$ and every pair of coprime numbers $(p,q)$, there is one simple geodesic loop of $(p,q)$ type through every vertex of $T$.
\end{abstract}


\section{Introduction}

There is a long history in the study of simple closed geodesics on compact Riemannian manifolds and also on non-regular surfaces, in particular on convex polyhedra, where may be the first systematic account is in the book of Alexandrov \cite{Al50}. There is a short review of this history in the introduction of \cite{BS20}. Recently, simple closed geodesics   
were studied by Protasov (see \cite{Pr07}) on tetrahedra in Euclidean space, Fuchs and Fuchs in \cite{FF07} proved that closed geodesics in regular tetrahedra have no self-intersections and gave nice characterizations of closed geodesics and the first author of this paper and Sukhorebska described simple closed geodesics on regular tetrahedra in the hyperbolic and spherical space (see \cite{BS19,BS20,BS21} and \cite{Bo22}). 

On a convex polyhedron (in Euclidean, spherical or hyperbolic space), a simple closed geodesic does not pass through any vertex of the polyhedra. Then, the study of simple loops which contain a vertex and are geodesics at every other point needs a separated study. These kind of loops are called simple geodesic loops (see the precise definitions in the next section). Their study for regular tetrahedra in the Euclidean space was done by Davids, Dods, Traub and Yang proved  in \cite{DDTY}. In fact, they proved that there are not geodesic loops on regular tetrahedra and cubes in the Euclidean Space $\re^3$.

In this paper we start the study of closed geodesic loops on tetrahedra in spherical and hyperbolic spaces. As it should be expected, the behaviour results different from the euclidean case. 

In fact, for the spherical space $\es^3$ we prove: 1) if the internal angles $\a_i$ of the faces of a tetrahedron (regular or not) satisfy $\pi/3 < \a_i<\pi/2$ then the tetrahedron has no simple geodesic loop (Th. \ref{th1S}); 2) for regular tetrahedra with $\a_i=\pi/2$  the result is still like in the Euclidean case: there are no simple geodesic loops on these tetrahedra; 3) For any tetrahedra (regular or not) with $\a_i\ge \pi/2$ and the lengths $a_i$ of the edges of the tetrahedron satisfying $a_i>\pi/2$, there are three simple geodesic loops through every vertex of the tetrahedron (Th. \ref{th2S}), a behavior that does not appear in the Euclidean case. 

For the hyperbolic space $\H^3$ we prove that for every regular tetrahedron $T$ and every pair of coprime numbers $(p,q)$, there is one simple geodesic loop of $(p,q)$ type through every vertex of $T$ (Th. \ref{th2S}). A result which is completely opposite to the corresponding Euclidean one.

A concept related to closed geodesics and closed geodesic loops is that of quasi-geodesic. A geodesic in a convex polyhedron is a geodesic on the ambient space on each face of the polyhedron, and at the points in the edge the angle between the two segments of geodesic starting from this point is $\pi$ from both sides  measured on the polyhedron. Then Alexandrov  defined quasi-geodesics on a polyhedron in the euclidean space $\re^3$ as curves which are lines on the faces of the polyhedron and  at the points of intersections with the edges the angles between the segments (measured on the polyhedron) are $\pi$ from both sides of the curve and at the points of intersection with vertices, the angles between the segments is lower or equal than $\pi$ from both sides of the curve (then, unlike geodesics, quasi-geodesics can pass through vertices). 

A. Pogorelov (\cite{Po49}) proved (by approximation of polyhedra by surfaces) that every convex polyhedron in the Euclidean space has at least three simple closed quasi-geodesics. This was an existence theorem which did not describe the quasi-geodesics. 

O'Rourke and Vilcu faced the problem of describing the quasi-geodesics on tetrahedra using purely geometric methods. In  \cite{ORVi22} they prove that every tetrahedron in the Euclidean space has a $2$-vertex quasi-geodesic, a $3$-vertex quasi-geodesic, and a simple closed geodesic or a $1$-vertex simple closed quasi-geodesic. Since, according our definitions, every $1$-vertex simple closed quasigeodesic is also a simple geodesic loop, and it was proved in \cite{DDTY} that there are no simple geodesic loops in regular tetrahedra, the possible $1$-vertex quasigeodesics only can be done in non regular tetrahedra.

The geodesic loops that we have found in the tetrahedra of spherical and hyperbolic spaces are also quasigeodesics. In this sense, our results also extend the study of $1$-vertex quasigeodesics for  tetrahedra to all simply connected spaces of constant sectional curvature.

In the next section we shall recall with more precision the definitions of the concepts we are working on and will state our main results. Along the other two sections we shall prove them.

\section{Definitions and main results}

\begin{defi}
A geodesic curve in a (non necessarily regular) surface embedded in a simply connected $3$-space of constant sectional curvature is a curve $\gamma$ that, for any point $x$ in its trace $\Gamma$, there is a neighborhood $U$ of $x$ in $\Gamma$ such that, for every two points $p,q$ in this neighborhood, the segment of $\gamma$ between $p$ and $q$ realizes the distance between $p$ and $q$ on the surface. This is an existence theorem, 
\end{defi}

 With this definition, it happens that on a convex polyhedron, a simple closed geodesic does not pass through any vertex of the polyhedra. In order to allow vertices, the following definition was introduced. 
 
 \begin{defi}
 Given a polyhedron, a simple (without self-intersections) closed curve through only one vertex $A$ which is a geodesic except at $A$ is called a simple geodesic loop.
 \end{defi} 
 
 As we said in the introduction, there are no geodesic loops on regular tetrahedra in the Euclidean space. For tetrahedra in the sphere $\es^3$ of constant sectional curvature $1$ we shall prove:
 
 \begin{teor}\lb{th1S}
Let $T$ be a tetrahedron in $\es^3$ with angles $\a_i$ at the vertices of its faces. If $\pi/3 < \a_i < \pi/2$ for every $\a_i$, then the tetrahedron  has no simple geodesic loop.
\end{teor}

\begin{teor}\lb{th2S}
Let $T$ be a tetrahedron in $\es^3$ with angles $\a_i$ at the vertices of its faces.
If $\pi/2 < \a_i$ for every angle $\a_i$ and $\pi/2< a$ for the length $a$ of every edge of the tetrahedron, then, for each vertex $A_i$ of the tetrahedron  there exist $3$ simple  geodesic loops through $A_i$. 
\end{teor}

For the study in the hyperbolic space we need also the following concepts.
 
 \begin{defi}
In a geodesic line in a tetrahedron, the points where the geodesic intersects the edges are called the vertices of the geodesic.

A simple closed geodesic (resp. a simple geodesic loop) in a tetrahedron in the hyperbolic space $\H^3$ is of type $(p,q)$ if it has $p$ vertices  on each of two opposite edges of the tetrahedron, $q$ vertices on each of other two opposite edges, and $p + q$ vertices on each of the remaining two opposite edges.
\end{defi}

Por tetrahedra in the hyperbolic space $\H^3$ of constant sectional curvature $-1$ we shall prove:

\begin{teor}\lb{th1H}
For every vertex $A$ of a regular tetrahedron $T$ in $\H^3$ and for every pair $(p,q)$ of coprime integers, there is a simple geodesic loop through $A$ which is of type $(p,q)$ and unique up to isometries.  
\end{teor}

{\it Moreover, every simple geodesic loop $\gamma$ is of type $(p,q)$ for some coprime pair $(p,q)$.}

\section{Simple geodesic loops on tetrahedra in $\es^3$}

In this section we shall prove theorems \ref{th1S} and \ref{th2S}. We shall consider non necessarily regular tetrahedra in the {\it sphere $\es^3$ of constant sectional curvature $1$}, and we denote the angles of the triangular faces of the tetrahedra by $\a_i$. 
First, we shall prove the following lemma:

\begin{lema}\lb{leS}
Let $T$ be a tetrahedron in $\es^3$ with vertices $A_1,A_2,A_3,A_4$. A simple closed geodesic loop through one vertex $A_i$ can only cut two edges of the tetrahedron, and only one time. These edges are two of the three that has the face opposite to $A_i$.
\end{lema}
\begin{demo} Let us suppose, without loose of generality, that $i=4$.

Let us consider the development of the faces of the tetrahedron on a sphere $\es^2$ of dimension $2$ and sectional curvature $1$, like in Figure \ref{fig1}. In this development, the points $A_{41}, A_{42}, A_{43}$ correspond to a unique point $A_4$ in the tetrahedron, and the  pairs of edges $A_iA_{4j}$ having the same value of $i$ are the same edge in the tetrahedron.

 A geodesic loop in the tetrahedron has to be the union of geodesic segments contained in its faces. Then, in the development, a closed geodesic loop at $A_4$ will start from one $A_{4j}$, let us say $A_{42}$. It will contain a geodesic segment $\gamma_2$ from $A_{42}$ to the edge $A_1A_3$, which will cut it at some point $P_2$. This segment cannot go to edges $A_1 A_{42}$ neither $A_3A_{42}$, because in the last two cases you will have two geodesic segments of length lower than $\pi$ joining two points in the sphere, which is impossible. 
 
 \begin{figure}[!h]
\begin{center}
\includegraphics[scale=0.7]{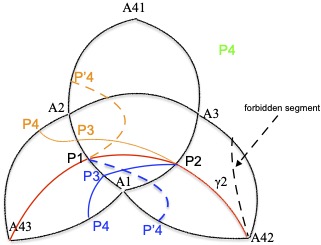}
\caption{}\lb{fig1}
 \end{center}
\end{figure}
 
 In the development, the geodesic loop have to finish in a geodesic segment with one end in $A_{43}$ or $A_{41}$. Let us suppose it is $A_{43}$, then the other end  of the geodesic segment must be a point $P_1$ of the edge $A_2A_1$. If the loop closes just by a segment joining $P_1$ and $P_2$, the statement of the  theorem is true. If not, there must be another segment of geodesic from $P_2$ to another point $P_3\ne P_1$ of the edge $A_1 A_2$. Since the loop cannot cross itself, it has to continue going to hit at some point $P_4$ in the unique edge among $A_1A_{43}$ and $A_2A_{43}$ which is possible without self-intersection. But this point is the same at the tetrahedron that a point $P'_4$ in the edge $A_1A_{42}$ or $A_2A_{43}$. But the union of the segments $A_{42}P_2$, $P_2 P_3$, $P_3P_4$ divides the development of the tetrahedron in two connected components in such a way that $P_4'$ is in a connected component and $P_1$ in the other, then they cannot join to close the loop without self-intersection.
 
The same argument can be used if the geodesic loop finishes with a geodesic segment with one end in $A_{41}$. Then the Lemma is proved.

\end{demo}

As a consequence of this Lemma, for the proof of  theorems \ref{th1S} and \ref{th2S}, we have to consider only simple closed geodesic loops through a vertex that cut two of the three edges of the face opposite to the vertex.

\ \\

\begin{demo}{\bf of Theorem \ref{th1S}}.

 Let us consider the development of the faces of the tetrahedron on a sphere $\es^2$ of dimension $2$ and sectional curvature $1$, like in Lemma \ref{leS}. Then the lengths of the geodesics from $A_{43}$ to $A_1$ and from $A_{42}$ to $A_1$ are the same. We shall write $|A_{43}A_1| = |A_{42}A_1|$.

On each one of the faces of the tetrahedron, with vertices $A_i$, angles $\a_i < \pi/2$ and opposite sides $a_i$, we can apply the second cosinus law 
\begin{align}
\cos \a_i = - \cos \a_j \cos \a_k + \sin\a_j\sin\a_k \cos a_i
\end{align}
to conclude that $a_i<\pi/2$, because all  other terms in the formula which do not contain $a_i$ are  positive, then $\cos a_i>0$. Then we conclude that $|A_{43}A_1| = |A_{42}A_1|< \pi/2$.
\begin{figure}[!h]
\begin{center}
\includegraphics[scale=0.7]{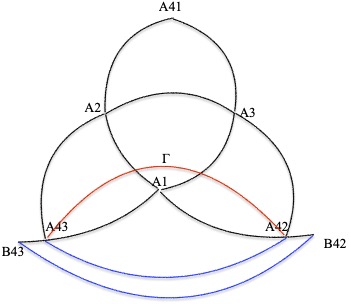}
\caption{}\lb{fig2}
 \end{center}
\end{figure}
Since $\pi/2 > \a_i > \pi/3$, then the angle $\pi/2 < \widehat{A_{43} A_1 A_{42}} < \pi$, then there is a geodesic $A_{43}A_{42}$ external to the development of the tetrahedron. We consider the triangle $A_{43} A_1 A_{42}$ which has two equal sides $|A_{43}A_1| = |A_{42}A_1| < \pi/2$, 
and $\widehat{A_1A_{43}A_{42}} = \widehat{A_1A_{42}A_{43}}$. Let us extend the geodesics $A_1 A_{43}$ and $A_1A_{42}$ to $A_1 B_{43}$ and $A_1B_{42}$, respectively, so that $|A_1 B_{43}| = |A_1B_{42}|= \pi/2$. The triangle $A_{43} A_1 A_{42}$ is contained in the triangle $B_{43} A_1 B_{42}$, and the angles $\widehat{A_1 B_{42}B_{43}} = \widehat{A_1 B_{43}B_{42}}=\pi/2$, then, by Gauss-Bonnet, $\widehat{A_1A_{43}A_{42}} = \widehat{A_1A_{42}A_{43}}< \pi/2$. By the law of sines, $|A_{43} A_{42}| < \pi$. If there is a geodesic loop at $A_4$, when we develop the tetrahedron it produces a geodesic $\Gamma$ from $A_{43}$ to $A_{42}$ cutting edges $A_1A_2$ and $A_1A_3$ (or similar situations with a geodesic from $A_{43}$ to $A_{41}$ or another from $A_{42}$ to $A_{41}$ ). Since $\widehat{A_1A_{43}A_{42}}<\pi/2$ and $\widehat{A_1A_{43}A_{2}}<\pi/2$, the angle at $A_{43}$ between the geodesic $\Gamma$ and the external one $A_{43}A_{42}$ is lower than $\pi$, and we have a similar bound at $A_{42}$. Then $\Gamma$ and $A_{43}A_{42}$ are two different geodesics cutting at two points at distance $|A_{43}A_{42}| <\pi$, which is impossible, then there is no geodesic loop at $A_4$.
\end{demo}

\begin{nota*} 
From Theorem \ref{th1S} it follows that regular tetrahedra with the interior angles of their faces $\a <\pi/2$ have no simple geodesic loop with $1$ vertex.

For the case of regular tetrahedra with $\a_i=\pi/2$, one has also $a_i=\pi/2$, and in the deveolpment done in the previous proof all the angles and sides are $\pi/2$. The points $A_{42}$ and $A_{43}$ are united by one external geodesic $A_{42}A_{43}$ and an internal geodesic formed for the union of the segments $A_{43}A_2$, $A_2A_3$ and $A_3A_{42}$, but this geodesic goes through $3$ vertices of the tetrahedron, then, also in this case, there are no geodesic loops containing only one vertex.
\end{nota*}

\begin{demo}{\bf of Theorem \ref{th2S}}

Let us call $A_4$ one of the vertices of the tetrahedron, and let us consider a development of the tetrahedron on $\es^2$ where the vertix $A_4$ developes in three vertices $A_{41}$, $A_{42}$ and $A_{43}$. Let us consider the bisectrix $b$ of the angle $\widehat{A_{42}A_1A_{43}}$ starting from $A_1$, and the point $O$ in the bisectrix such that $|OA_{42}|=\pi/2 = |OA_{43}|$. From the cosinus law applied to the triangle $A_1OA_{43}$, we have
\bec\lb{cosl1}
0= \cos|OA_{43}| = \cos|A_{43}A_1| \cos|OA_1| + \sin|A_{43}A_1| \sin|OA_1| \cos\widehat{OA_1A_{43}}.
\eec
Since the interior angles of the triangle are bigger than $\pi/2$, then $\widehat{OA_1A_{43}} < \pi/4$. Moreover $|A_1A_{43}| >\pi/2$ by hypothesis. Then the second addend in the last term of the equality \eqref{cosl1} is positive, and $\cos|A_{43}A_1|<0$, which implies $\cos|OA_1|>0$, then  $|OA_1| <\pi/2$. 
\begin{figure}[!h]
\begin{center}
\includegraphics[scale=0.5]{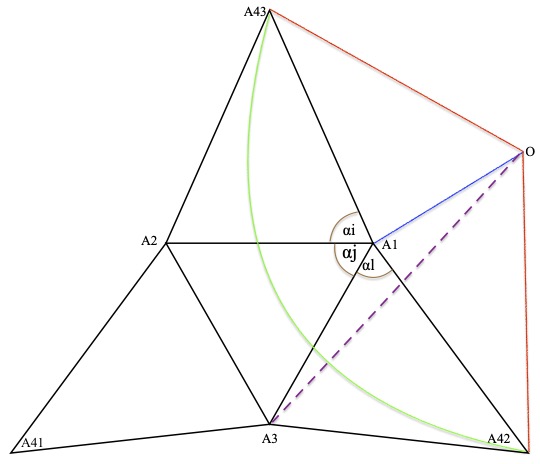}
\caption{}\lb{fig3}
 \end{center}
\end{figure}

Now, we consider the segment of circle $C$ from $A_{42}$ to $A_{43}$ and center at $O$. Because its radius is $\pi/2$, it is a segment of geodesic. If we prove that $|OA_3| > \pi/2$ and $|OA_2| > \pi/2$, then $C$ will be included in the development of the tetrahedron and a $C$ wiil become a loop when we reconstruct the tetrahedron from its development. Let us prove $|OA_3| > \pi/2$. We consider now the triangle $OA_1A_3$. The cosinus law for this triangle gives
\bec\lb{cosl2}
\cos|OA_{3}| = \cos|A_{3}A_1| \cos|OA_1| + \sin|A_{3}A_1| \sin|OA_1| \cos\widehat{OA_1A_{3}},
\eec
but $\widehat{OA_1A_{3}} > \pi/2$, $|OA_1| < \pi/2$ and $|A_{3}A_1| >\pi/2$, then the second term in \eqref{cosl2} is negative, which implies $|OA_3|>\pi/2$. The same argument proves  that $|OA_2| > \pi/2$.

 This geodesic $C$ is the unique inside the development joining $A_{42}$ to $A_{43}$. In fact, 
 it follows from the cosinus law applied to the triangle $A_{42}A_1A_{43}$ that $|A_{42}A_{43}|<\pi/2$, then the unique geodesics in $\es^2$ joining  $A_{43}$ and $A_{42}$ are $C$ and the geodesic $D$ which completes $C$ to give a closed geodesic of $\es^2$, but $D$ is outside of the development of the tetrahedron, then $C$ is unique. 

We can obtain other two geodesic loops through $A_4$ by repeating the construction for the angles $\widehat{A_{43}A_2A_{41}}$ and $\widehat{A_{41}A_3A_{42}}$. This finishes the proof of Theorem \ref{th2S}.
\end{demo}

\begin{nota*}
Regular tetrahedra  with the interior angles of their faces $\a>\pi/2$ also satisfy that $a>\pi/2$, then they satisfy the conditions of  Theorem \ref{th2S} and they have three simple geodesic loops through each vertex.
\end{nota*}

\section{Simple geodesic loops on tetrahedra in $\H^3$}

In this section we shall consider non necessarily regular tetrahedra in the {\it sphere $\H^3$ of constant sectional curvature $-1$}.

\begin{lema}\lb{dbound}
Let $T$ be a  tetrahedron in $\H^3$ such that for every three faces with a common vertex $V$, the sum $2 \widetilde \a$ of the angles of the faces at $V$ is lower than $\pi$.  Let $\gamma$ be a geodesic except may be  at some point $A_0$ in an edge, such that there is a point $Q\in \gamma-\{A_0\}$ realizing the distance $d$ between a vertex $A\neq A_0$ of $T$ and $\gamma$. If  $\sinh d \le \cos\widetilde\a \sin h$, where $h$ is the minimum of the distances between $A$ and its opposite edges in $T$, then  $\gamma$ has self-intersections.
\end{lema}
\begin{demo}
Let $T$ be a  tetrahedron with edges of length $\ell_i$ and angles at the vertices of their faces equal to $\a_i$ (remember $0<2 \widetilde \a = \sum \a_i<\pi$) and height $h$.   

Given a vertex $A$ of $T$, let $h$ the minimum distance from $A$ to their opposite edges. We consider the union $U$ of the three faces of $T$ with common vertex $A$, a geodesic $g$ of one of these faces starting from $A$, and the development $D$ of $U$ on $\H^2$ cutting $U$ along $g$. In the boundary of $D$, there are two segments of geodesic $g_1, g_2$ of length $|g_1| = |g_2|=|g|\ge h$, where the inequality follows from the definition of $h$. Obviously, the sector of a disk between $g_1$ and $g_2$ with center at $A$ and radius $h$  is contained in $D$. Let $P_1 \in g_1$ and  $P_2 \in g_2$ points at distance $h$ from $A$. Since the geodesic disks in $\H^2$ are convex, the geodesic $\gamma_{12}$ joining $P_1$ and $P_2$ is contained in $D$. 
\begin{figure}[H]
\begin{center}
\includegraphics[scale=0.7]{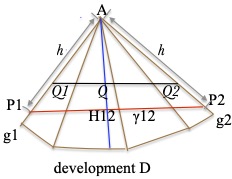}   
\caption{}\lb{fig4}
\end{center}
\end{figure}
Moreover, the geodesic bisectrix $\delta_{12}$ of the angle $\widehat{P_1AP_2}$ cuts orthogonaly $\gamma_{12}$ at some point $H_{12}$ and realizes the distance between $A$ and $\gamma_{12}$. Let us observe also that $\widehat{P_1AP_2} = \sum \a_i<\pi$. We now consider the triangle ${H_{12}AP_2}$, with sides of lengths $h$, $b_{12}:=|P_1 P_2|/2$ and $d_{12} =|AH_{12}|$ and angles $\widetilde\a:= \sum \a_i/2<\pi/2$ at $A$ and $\pi/2$ at $H_{12}$. By the sinus and cosinus laws applied to this triangle, we obtain:

 \begin{align}
 {\sinh(b_{12})} & = \sinh(h) {\sin \widetilde\a}, \quad
 \cosh h = \cosh d_{12} \cosh(b_{12})
 \end{align}
 From these to equalities it follows that
 \begin{align}
 \sinh(h) &= \frac{\sinh d_{12}} {\cos(\widetilde\a)}\quad  \text{ and } \quad 
 \sinh b_{12}=  {\sinh d_{12}}\ {\tan(\widetilde\a)}\lb{reldlb}
 \end{align} 
Let $Q_{12}$ be a point in $\delta_{12}$ such that $|AQ_{12}|<d_{12}$, and let $\overline \gamma$ be a geodesic of $\H^2$ orthogonal to $\delta_{12}$. Then $|AQ_{12}|$ is the distance from $A$ to $\overline\gamma$, and  $\overline\gamma$ will contain the segment of geodesic between $Q_1\in g_1$ and $Q_2\in g_2$ such that $|AQ_1| = |AQ_2|$.

 Let us consider the loop $\gamma$ in the statement of the theorem, let $\delta$ be the geodesic starting from $A$ that realices the distance between $A$ and $\gamma$, and let $g$ be the geodesic starting from $A$ such that when we do the development $D$ of $U$ cutting $U$ along $g$, the geodesic $\delta$ is the bisectrix of the angle of $D$ at $A$. 
 
 In the development $D$, the  loop $\gamma$ is the union of geodesic segments in $\H^2$ with ends in the boundary of $D$.  Let $\xi$ be the segment of $\gamma$ where is the point $Q\in\gamma$  which realizes the distance $d$ from $\gamma$ to $A$.
 
 The condition  $\sinh d \le \cos\widetilde\a \sinh h = \sinh d_{12}$ implies (by the uniqueness o the geodesic through one point in a given direction) that this segment must be part of a geodesic $\overline \gamma$ like the  one considered in the previous paragraph, and its ends $Q_1$ and $Q_2$ are at the same distance from $A$, then they are the same point in $T$, and the loop $\gamma$ has a self-intersection at $Q_1\sim Q_2$ in $T$. 
\end{demo}


\begin{demo}{of Theorem \ref{th1H}}

Let us consider an edge $A_1A_2$ of the tetrahedron. From \cite{BS20} we know that there is a point  $A_0$ in it, at some distance $a$ from $A_1$ such that  there is a simple closed geodesic of type $(p,q)$ from $A_0$ to $A_0$. Let $D(p,q)$ be the development of the tetrahedron obtained by unrolling this geodesic on the hyperbolic plane $\H^2$. In $D(p,q)$ $A_1$ is represented by two points $A'_1$ and $A"_1$. 
\begin{figure}[H]
\begin{center}
\includegraphics[scale=0.6]{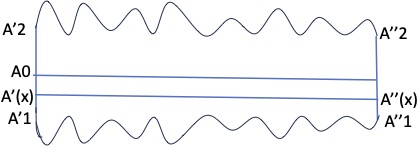}   
\caption{}\lb{fig5}
\end{center}
\end{figure}
 Now, let us move the point $A_0$ to a point $A(x)$ on the edge $A_1A_2$ at distance $x$ in the direction of $A_1$ ($0<x<a$). In $D(p,q)$, each $A(x)$ is represented by two points $A'(x)$ and $A"(x)$. For small values of $x$, there is a geodesic $\gamma_x$ in $\H^2$, contained in $D(p,q)$ joining $A'(x)$ with $A"(x)$. In the tetrahedron $T$ this curve gives a simple geodesic edge loop (a loop which is geodesic except may be at a point $A(x)$ in the edge $A_1A_2$). Let $x_0\le a$ be the supremum of the $x$ for which the above $\gamma_x$ exists and is contained inside the development. 
 
 If $x_0=a$, then $\gamma_a$ is the simple geodesic loop at $A_1$ that we are looking for. 
 
 Let us suppose that $x_0<a$. Then $\gamma_{x_0}$ will hit some vertex $A$ in $D(p,q)$ different from $A'_1$ and $A"_1$. For $\eta>0$ small enough, there is an $\eps$ ($0<\eps<a$) such that $\dist(A,\gamma_{x_0-\eps}) = \eta$. By definition of $x_0$, $\gamma_{x_0-\eps}$ has no self-intersection.
 
For every vertex $V\in D(p,q)-\{A_1\}$, and $x\in [0,x_0]$, let $c(V,x)= |A(x)V|$. The function $x\mapsto c(V,x)$ is continuous in $[0,x_0]$ and bounded from below  by a positive number $c(V)$. Since the number of vertices $V$ is finite, there is an infimum $c>0$ among all the $c(V)$.

By definition of $c$, $\dist(A,A(x_0-\eps) > c$. Let $Q$ be the point in $\gamma_{x_0-\eps}$ realizing the distance $\eta$ between  $A$ and $\gamma_{x_0-\eps}$. If we take $\eta < c$ then $Q\ne A(x_0-\eps)$. Taking $A(x_0-\eps)$ as the $A_0$ and $\gamma_{x_0-\eps}$ as the $\gamma$ in Lemma \ref{dbound}, the general conditions of that lemma are satisfied. Moreover, if we choose $\eta$ such that it also satisfies $\sinh \eta <\cos\widetilde\a \sinh\zeta$, with $\zeta = \min\{h,c\}$, the condition on $\sinh d$ of the Lemma \ref{dbound}  is also satisfied. Then, by that lemma, $\gamma_{x_0-\eps}$ have self-intersections, in contradiction with our previous conclusion. 

Then $x_0=a$, and $\gamma_a$ is a simple geodesic loop at vertex $A_1$, and only contains this vertex.

By construction, $\gamma_a$ cuts the same sides that $\gamma$ in the development $D(p,q)$, then it is of type $(p,q)$.

Since the simple closed geodesic used to construct the geodesic loops are unique up to isometries in the tetrahedron, the geodesic loops that we have constructed are also unique up to isometries of the tetrahedron. 
\end{demo}

Moreover, the simple geodesic loops constructed before are all the possible simple geodesic loops.  In fact, if $\gamma$ is a geodesic loop at some vertex $A_1$, it is the union of segments of geodesics on the facets of the tetrahedron, and Protasov proved in \cite{Pr07} that all simple closed curves of this form are of type $(p,q)$ for some pair of coprime numbers. Let $D(p,q)$ the development of the tetrahedron obtained unrolling $\gamma$. Let $A_2$ the other vertex of the edge $A_1A_2$ which appears as the first and the last in the development $D(p,q)$. We know, by \cite{BS20} that there is a simple closed geodesic starting at some point of the edge $A_1A_2$ giving the same development $D(p,q)$. Application of our previous construction to this geodesic will give a loop $\gamma_1$ at $A_1$ which, in the development $D(p,q)$, is a geodesic in $\H^2$ joining the points $A'_1$ and $A"_1$ in $D(p,q)$. By the uniqueness of a geodesic joining two points in the hyperbolic plane, $\gamma_1=\gamma$. \\

If $T$ is a non necessarily regular tetrahedron in $\H^3$ with angles of their faces $\a_i\le \pi/4$, the development $D(p,q)$ of a simple closed geodesic is a convex polygon in $\H^2$ (see \cite{BS20}, section 6), then, in the development at the beginning of the previous proof, the existence of the geodesic from $A'_1$ to $A"_1$ is granted by this convexity, and we have:

\begin{prop}
\lb{pr1H}
For every vertex $A$ of a tetrahedron $T$ in $\H^3$ with angles of their faces $\a_i\le \pi/4$ and for every pair $(p,q)$ of coprime integers, there is a simple geodesic loop through $A$ which is of type $(p,q)$ and unique up to isometries.
\end{prop}

\bibliographystyle{alpha}

\begin{thebibliography}{99}
 
 \bibitem{Al50} A. D. Alexandrov, {\it Convex polyhedra}, Moscow–Leningrad, GITTL 1950, 428 pp.;  English transl. in Springer Monogr. Math., Springer-Verlag, Berlin 2005, xii+539 pp.

\bibitem {Bo22}   Borisenko, A. A. A necessary and sufficient condition for the existence of simple closed geodesics on regular tetrahedra in spherical space. {\it Mat. Sb.} 213 (2022), no. 2, 37–49; English transl. in {\it Sb. Math. }213 (2022), no. 2, 161–172 

\bibitem{BS19} Borisenko, A. A.; Sukhorebska, D. D. Classification of simple closed geodesics on regular tetrahedra in a Lobachevskiĭ space. (Russian) {\it Dopov. Nats. Akad. Nauk Ukr. Mat. Prirodozn. Tekh. Nauki} 2019, no. 4, 3–9. 

\bibitem{BS20} Borisenko, A. A.; Sukhorebska, D. D. Simple closed geodesics on regular tetrahedra in  Lobachevskiĭ space. {\it Mat. Sb.} 211 (2020), no. 5, 3–30; English transl. in
Sb. Math. 211 (2020), no. 5, 617–642.


\bibitem{BS21} Borisenko, A. A.; Sukhorebska, D. D. Simple closed geodesics on regular tetrahedra in spherical space. {\it Mat. Sb.} 212 (2021), no. 8, 3–32; English transl. in {\it Sb. Math. }212 (2021), no. 8, 1040–1067 

\bibitem{DDTY}  Davis, Diana; Dods, Victor; Traub, Cynthia; Yang, Jed Geodesics on the regular tetrahedron and the cube. {\it Discrete Math.} 340 (2017), no. 1, 3183–3196

\bibitem{FF07}  Fuchs, Dmitry; Fuchs, Ekaterina, Closed geodesics on regular polyhedra. {\it Mosc. Math. J.} 7 (2007), no. 2, 265–279, 350.

\bibitem{ORVi22} Joseph O’Rourke  and Costin Vîlcu, Simple Closed Quasigeodesics on Tetrahedra. {\it Information} 2022, 13, 238, page 1 to 20.

\bibitem{Po49} Pogorelov, A.V. Quasi-geodesic lines on a convex surface. {\it Mat. Sb.} 1949, 25, 275–306; English transl. in {\it Am. Math. Soc. Transl.} 1952, 74.

\bibitem{Pr07} V. Yu. Protasov, Closed geodesics on the surface of a simplex, {\it Mat. Sb.} 198:2 (2007), 103–120; English transl. in {\it Sb. Math.} 198:2 (2007), 243–260.

\end{thebibliography}

{ Alexander Borisenko, \\
{B. Verkin Institute for Low Temperature, Physics and Engineering of the National Academy of Sciences of Ukraine\\
Kharkiv, Ukraine
}\\ and \\
 Department of Mathematics \\
University of Valencia \\
46100-Burjassot (Valencia), Spain}

{oleksandr.borysenk@uv.es}\\

{ Vicente Miquel \\
Department of Mathematics \\
University of Valencia \\
46100-Burjassot (Valencia), Spain
}

{vicente.f.miquel@uv.es}




\end{document}